\documentstyle[12pt]{article}
\makeatletter
\renewcommand{\@cite}[2]{[{{\bf #1}\if@tempswa,#2\fi}]}
\renewcommand{\@biblabel}[1]{[#1]\hfill}
\makeatother
\newtheorem{defin}{Definition}
\newtheorem{prop}{Proposition}

\newtheorem{th}{Theorem}
\newtheorem{lemma}{Lemma}

\newfont{\sdbl}{msbm9}
\newfont{\dbl}{msbm10 at 12pt}
\newcommand{\eqdef}{\stackrel{\rm def}{=}}
\newcommand{\proof}{{\bf Proof\ }}

\newcommand{\oo}{{\cal O}}
\newcommand{\ff}{{\cal F}}
\newcommand{\g}{{\cal G}}

\newcommand{\Lim}{\mathop {\rm lim}}

\newcommand{\Spec}{\mathop {\rm Spec}}
\newcommand{\Proj}{\mathop {\rm Proj}}
\newcommand{\Frac}{\mathop {\rm Frac}}
\newcommand{\dm}{\mathop {\rm dim}}
\newcommand{\rank}{\mathop {\rm rank}}

\newcommand{\dz}{{\mbox{\dbl Z}}}
\newcommand{\sdz}{{\mbox{\sdbl Z}}}
\newcommand{\dbp}{{\mbox{\dbl P}}}

\newcommand{\Image}{{\rm Im}\:}

\newcommand{\f}{{\cal F}}

\begin{document}
\author{D. V. Osipov}
\title{Stability of torsion free sheaves on curves and infinite-dimensional
Grassmanian manifold}
\date{}
\maketitle

In~\cite{M} and~\cite{PW}
are described the Krichever map from the set of torsion free
sheaves on curves with some additional datas to
an infinite-dimensional Grassmanian manifold.
In this note we investigate the images of (semi)stable
rank 2 torsion free sheaves on curves via the Krichever map
and verify some analog of G.I.T. Hilbert-Mumford numerical
criterion (theorem~\ref{th1}) with respect to  actions of some one-parametric
subgroups of the group $SL(2,k[[z]]$ on the determinant bundle
of the infinite-dimensional Grassmanian manifold.

\section{Krichever map}
Consider the following geometric datas:
\begin{enumerate}
\item $C$ is a reduced irreducible complete algebraic curve defined over
a field $k$.
\item $p$ is a smooth $k$-rational point.
\item $\ff$ is a torsion free coherent sheaf of $\oo_C$-modules on $C$ of rank
$2$.
\item $t_p$ is a local parameter of point $p$, i.~e.
$\hat{\oo_p} = k[[t_p]]$, where $\hat{\oo_p}$ is the completion of local
ring $\oo_p$ along the maximal ideal.
\item $e_p$ ia a basis of  rank $2$ free  module $\hat{\ff_p}$
over the ring  $\hat{\oo_p}$, where $\hat{\ff_p} \eqdef \ff \otimes_{\oo_C}
\hat{\oo_p} $.
\end{enumerate}

In the sequel we will {\em call} such collection $(C,p,\ff,t_p,e_p)$
a {\em quintet}. If the Euler characteristic $\chi(\ff)= \mu$,
then we will say, that the quintet have an index $\mu$.

On the other hand, consider $V= k((z)) \oplus k((z))$,
$V_0 = k[[z]] \oplus k[[z]]$. Then {\em define} the infinite Grassmanian
$Gr^{\mu}(V)$ as the set of $k$-vector subspaces $W$ of $V$,
which are Fredholm of index $\mu$, i.~e.,
$\dm_k V / (V_0 + W ) < \infty $, $\dm_k W \cap V_0 < \infty$, and \\
$\dm_k W \cap V_0   - \dm_k V / (V_0 + W )  = \mu   $.

There exists {\em the Krichever map} $K$ from the set of quintets of index
$\mu$ to $Gr^{\mu}(V)$, which can be shortly defined as the following
chain of
evident maps:
$$
(C,p,\ff,t_p, e_p) \to
 H^0(C \backslash p, \ff)
\hookrightarrow
H^0 ( \Spec \oo_p
\backslash p,\ff)
\hookrightarrow
H^0 ( \Spec \hat{\oo_p}
\backslash p,\ff) =  \qquad \qquad
$$
\begin{equation}         \label{chain}
\qquad \qquad \qquad  \qquad      \qquad
=
 \hat{\ff_p} \otimes_{\hat{\oo_p}} \hat{K_p}
\stackrel{e_p}{\to}
(\hat{\oo_p} \oplus \hat{\oo_p})
\otimes_{\hat{\oo_p}} \hat{K_p}
\stackrel{t_p}{\to}  k((z)) \oplus k((z))  \; \mbox{,}
\end{equation}
where $\hat{K_p}$ is the fraction field of  $\hat{\oo_p}$.

The Krichever map $K$ and its basic properties is considered in details by
many authors. See, for example, \cite{M} for an algebraic description.
(But note that
the definition in~\cite{M} is slightly different from our definition.)
And see~\cite{PW} for an analytic description.

For any $W \in Gr^{\mu}(V) $  define the ring
$$A_W  \eqdef \{f \in k((z)) : fW \subset W \} $$

\begin{lemma}[see~\cite{M}]    \label{odin}
For any $W \in Gr^{\mu}(V) $
and any subring $A \subset A_W$  such that $k \subset A$   we have

1) $A \cap k[[z]] = k$   and

2) if $A \ne k$, then the ring $A$ has dimension $1$ over $k$.
\end{lemma}
\proof.
The first statement follows from $\dm_k W \cap V_0 < \infty$.

For the proof of the second statement consider the following subgroup of
$\dz$:
$$
\left\{ n \in \dz  \; \, : \; \, n= \nu(a) \, , \; a \in  \Frac A
 \right\} \mbox{,}
$$
where $\nu$ is the discrete valuation of $k((z))$
and $\Frac A$ is the fraction field of the ring $A$.
Then this subgroup is $r \dz$ for some integer  $r$.
Therefore

1) there exist $f, g \in A$ such that $r = \nu(f) - \nu(g)$

2) for any integer $n$
\begin{equation}  \label{zvezda}
\dm\nolimits_k A \cap z^{-nr}k[[z]] / A \cap  z^{(-n+1)r}k[[z]] \le 1
\end{equation}
Now from~(\ref{zvezda}), the first statement of this lemma,
and $r = G.C.D. (\nu(f), \nu(g))$
we have
$$
\dm\nolimits_k A / k[f,g] < \infty    \quad \mbox{.}
$$
Thus the dimension of the ring $A$ over $k$ is at most two.
Now let us assume that $f$ and $g$ are algebraically independent
elements over $k$. Then   for all pairs of integers $n_1 > 0$
and $n_2 > 0$ elements $f^{n_1} g^{n_2}$
are included in one basis of vector space $A$ over $k$.
But it contradicts to~(\ref{zvezda}) and this contradiction
 concludes the proof of  lemma. \vspace{5pt}

Notice also the following well-known statement.
\begin{lemma}  \label{Shurik}
$W \in Gr^{\mu}(V) $
is in the image of the Krichever map
if and only if $\rank A_W =1$, i.e.,
there exist two elements $f, g \in A_W$ such that $\nu(f)$
and $\nu(g)$ are relatively prime.
\end{lemma}
\proof.
Note that the condition $\rank A_W = 1$ is equivalent to
\begin{equation}  \label{zvzvezda}
\dm\nolimits_k k((z)) / (A_W + k[[z]]) < \infty  \quad \mbox{.}
\end{equation}

Let  $W= K(C,p,\f,t_p,e_p)$. By $A_C$ denote the image of the ring
$H^0 (C \backslash  p, \oo_C)$  in $k((z))$:
$$
H^0 (C \backslash  p, \oo_C)
\hookrightarrow
\hat{K_p}
\stackrel{t_p}{\to}
k((z))    \quad \mbox{.}
$$
From the  Riemann-Roch theorem we have
$$
\dm\nolimits_k k((z)) / (A_C + k[[z]])   < \infty   \quad \mbox{.}
$$
From $A_C \subset A_W$  it follows~(\ref{zvzvezda}).

Now let $\rank A_W = 1$.
Fix any subring $A \subset A_W$ such that  $k \subset A$
and $\dm\nolimits_k A_W / A < \infty$.
By lemma~\ref{odin} the ring $A$ has dimension $1$ over $k$.
Then from graded $k$-algebra
$$
A_* \eqdef \bigoplus_{i=0}^{\infty} A \cap z^{-i}k[[z]]
$$
we construct the complete irreducible reduced curve
$C = \Proj A_*$.
And from~(\ref{zvzvezda}) we see
that the graded $A_*$-module
$$
W_* \eqdef
\bigoplus_{i=-\infty}^{\infty} W \cap z^{-i}V_0
$$
determines the rank 2 torsion free coherent sheaf $W_* \,\tilde{}$ on $C$.
Moreover, we have $C = \Spec A \cup p$,
where $p$ is a smooth $k$-rational point, $\hat{\oo_p}= k[[z]]$,
$\ff \mid_{C \backslash p} = W \, \tilde{}$,
$\hat{\ff_p} = V_0 $,
$t_p = z $,
$e_p = {(1,0), (0,1) \in V}$.
And
$$
\bigoplus_{i=0}^{\infty} H^0(C, \oo(ip)) \simeq A_*
\qquad
\mbox{,}
\qquad
 \bigoplus_{i=-\infty}^{\infty} H^0(C, \ff(ip)) \simeq W_*  \quad \mbox{,
and}
$$
$$
H^0(C, \ff) \simeq W \cap V_0
\qquad
\mbox{,}
\qquad
H^1(C, \ff) \simeq V / (W + V_0)  \quad \mbox{.}
$$
This lemma is proved.  \vspace{7pt}

\noindent {\bf Remark.} If we fix a triplet $(C, p, t_p)$
and identify $(C,p,\ff, t_p, e_p)$
with $(C, p, \ff', t_p, e_p')$
for every sheaf isomorphism $\alpha : \ff \to \ff'$
such that $\alpha (e_p) = e_p'$;
then the Krichever map is an injective map.
In this case $W \in  Gr^{\mu}(V)$ is in the image of such map if and only
if $A_C W \subset W$.  \\[4pt]

Now let $W =
K(C,p,\f,t_p,e_p) \in
   Gr^{\mu}(V)   $.
 Then with every subsheaf $\g$  of $\ff$
 we can associate  the $k((z))$-vector subspace
$L_{\g} \subset V$ such that $L_{\g} \cap W \ne 0$ :
$$
H^0(C \backslash p, \g) \hookrightarrow H^0(C \backslash p, \ff)
= W  \subset Gr^{\mu}(V)
$$
\begin{equation}  \label{zv}
\qquad
\g  \longmapsto  H^0(C \backslash p, \g) \cdot  k((z))  \subset{V}
\end{equation}
We have the following  proposition.

\begin{prop}    \label{p1}
Let  $W= K(C,p,\f,t_p,e_p) $,
$
R(\f)$ be the set of all rank 1 coherent subsheaves $ \g \subset \ff$
 such that the sheaf $ \ff/\g $ is a free torsion coherent
 sheaf,
$
R(W) $ be the set of all 1-dimensional $k((z))$-vector subspaces $L \subset
V$
such that $L \cap W \ne 0$.
Then  map~(\ref{zv})  is an one-to-one correspondence between
$R(\ff)$ and $R(W)$.\\
 Moreover, if $\g \in R(\ff) \mapsto L_{\g} \in R(W)$,
then
\begin{equation}  \label{ech}
\chi(C, \g) = \dm\nolimits_k W \cap L_{\g} \cap V_0
- \dm\nolimits_k L_{\g}/ (W \cap L_{\g} + L_{\g}
\cap V_0 )
\quad \mbox{.}
\end{equation}
\end{prop}
\proof.

It is not difficult to see that
the set $R(\f)$ is in an one-to-one correspondence with $1$-dimensional
$k(\eta)$-vector subspaces  of $\ff_{\eta} \eqdef H^0(\Spec k(\eta),\ff) $,
where  $\eta \hookrightarrow C$ is the general point.
We have the canonical imbedding of $\ff_{\eta}$
into $V$ as the third term $H^0 (\Spec \oo_p \backslash p, \ff) = \ff_{\eta}$
of chain~(\ref{chain}). And this imbedding gives us a one-to-one
correspondence between $R(\ff)$ and $R (W)$.

Moreover, we have the following explicit construction. Let $L \in R(W)$.
Then the graded $(A_C)_*$-module $(L \cap W)_* \eqdef \bigoplus\limits_{i =
-\infty}^{\infty} (L \cap W) \cap z^{-i} V_0$
determines the rank $1$
torsion free sheaf $\g_L = (L \cap W)_* \, \tilde{}$
on $C = \Proj (A_C)_*$. Then $L \to \g_L$
inverts
map~(\ref{zv})
and
$$H^0(C, \g_L) \simeq W \cap L \cap V_0  \quad \mbox{,}$$
$$H^1 (C, \g_L)  \simeq L / (W \cap L + L \cap V_0) \quad \mbox{.}$$

\section{Determinant bundle}
As was explicitly described in~\cite{AMP},
$  Gr^{\mu}(V)   $  admits a structure of an  algebraic scheme,
which represents a functor, which generalizes the usual finite-dimensional
grassmanian. There exists an open covering of $  Gr^{\mu}(V)   $
by means of  $  Gr^{\mu}_A(V)   $, where $A$ is a commensurable with $V_0$
$k$-vector subspaces in $V$, i.~e.,  $\dm_k (A+V_0)/ (A \cap V_0)  < \infty
$
and
$$  Gr^{\mu}_A(V)  \eqdef \{  W \in \  Gr^{\mu}(V)  :
V = A \oplus W \}
$$
is an infinite-dimensional affine space which is isomorphic
to the spectrum of a polynomial ring of infinitely many variables.
One can define a linear bundle (determinant bundle)  $Det$
on $Gr^{\mu}(V)$ such that for any  $W \in Gr^{\mu}(V)$
$$
Det_W  = \bigwedge^{max} (W \cap V_0) \otimes \bigwedge^{max} (V / W + V_0)^*
\qquad \mbox{.}
$$
Now let us describe the determinant bundle $Det$ more explicitly
("in coordinats"). For this goal we
identify the $k$-vector spaces $V$ and $k((t))$
by means of the following continuous isomorphism:
$$
V \longrightarrow k((t))
$$
\begin{equation}  \label{diez}
(\sum a_i z^i )  \oplus (\sum b_i z^i)   \longmapsto  \sum_i (a_i t^{2i}
+ b_i t^{2i+1})
\end{equation}
Note that $V_0 \to k[[t]]$.

Let $S_{\mu}$ be the set of sequences $\{ s_{-\mu +1}, s_{-\mu +2},
\ldots  \}$  of integers such that:

1) these sequences are
strictly decreasing;

2) $s_k = -k$  for $k \gg 0$.

Note that for every $\mu \in \dz  $  the set $S_{\mu}$
is in an one-to-one correspondence with the Young diagramms.

For every $S = \{ s_{-\mu +1}, s_{-\mu +2},
\ldots  \} \in S_{\mu}$  let
$$H_S  \in Gr^{\mu}(V)  $$
be the  $k$-vector space generated by $\{ t^{s_k}   \}$,
$$A_S  \subset k((t))$$
be the smallest  closed $k$-vector  space generated by
$\{  t^q \mid q \in \dz \backslash S               \} $.
Then
\begin{equation}   \label{dbldiez}
Gr^{\mu}(V)  = \bigcup_{S \in S_{\mu}} Gr_{A_S}^{\mu}(V)  \quad \mbox{.}
\end{equation}

For any $S \in S_{\mu}$ let $p_S : k((t)) \to H_S$
be the projection from the decomposition $ H_S \oplus A_S  = k((t))$.
By $0(\mu) \in S_{\mu}$ denote $\{ \mu-1, \: \mu -2,  \: \mu -3, \: \ldots \}$.

We are now going to describe the determinant bundle $Det$
in other words (see~\cite[~\S 7.7]{PS}). Given $W \in Gr^{\mu}(V)$
let us  define an {\em admissible} isomorphism to be an isomorphism
$$
w:   H_{0(\mu)} \longrightarrow W   \qquad \qquad   (w(t^{\mu-1}) \in W,
\;
w(t^{\mu-2}) \in W, \ldots  )
$$
such that
 $$p_{0(\mu)} \cdot w - Id  \quad :  \quad  H_{0(\mu)} \to H_{0(\mu)}$$
has finite rank. (Note that for any $W \in Gr^{\mu}(V)$ we can always
find the corresponding admissible isomorphism.)
Then the determinant bundle $Det$ consists of
$$ \{ [\lambda, w] : \lambda \in k,\; w \; \mbox{is an admissible isomorphism}
  \}   \quad \mbox{,}$$
which is identified with respect to the following equivalence relation:
$$
[\lambda, w]  \sim [\lambda \cdot \det (w'^{-1} \cdot w), w'] \quad \mbox{,}
$$
where $w'$ and $w$ are two admissible isomorphisms for $W \in Gr^{\mu}(V)$.
(Note that
$$w'^{-1} \cdot w - Id \quad : \quad  H_{0(\mu)} \to H_{0(\mu)}  $$
has  finite rank as an endomorphism of $H_{0(\mu)}$, and
hence $\det (w'^{-1} \cdot w)$ is well defined with respect to the basis:
$\{ t^{\mu -1}, t^{\mu-2}, t^{\mu-3}, \ldots \}$.)

Let $Det^*$ be the dual vector bundle to $Det$. For every
$S \in S_{\mu}$ we can  construct the global section
$$
\pi_S \in H^0(Gr^{\mu}(V), Det^*)  \qquad : \quad Det \longrightarrow k
$$
$$
[\lambda,w] \stackrel{\pi_S}{\longrightarrow} \lambda \det(p_S \cdot w)
\in k    \quad \mbox{,}
$$
where   $ p_S \cdot w - Id \; : \; H_0 \to H_S$  has  finite rank.
Therefore we can compute $\det (p_S \cdot w)$ with respect to
basises
$$
 \{t^{\mu -1}, \; t^{\mu -2}, \; \ldots    \} \in H_{0(\mu)}
\qquad \mbox{and} \qquad  \{t^{s_{-\mu +1}},\; t^{s_{-\mu+2}},\; \ldots   \}
\in H_S   \quad \mbox{.}
$$
And as proven in~\cite{AMP} and~\cite{PS},
we have a closed imbedding ("Pl\"ucker imbedding"):
$$
Gr^{\mu}(V) \hookrightarrow   \dbp (\Pi (S_{\mu})^*)  \eqdef \Proj (Sym(\Pi
(S_{\mu})))
$$
$$
W \longmapsto \left\{ \pi_S \to \pi_S (w) \eqdef \det (p_S \cdot w)  \right\}
\quad \mbox{,}
$$
where $\Pi (S_{\mu}) \subset H^0(Gr^{\mu}(V), Det^*)$
is the $k$-vector subspace generated by the global sections
$\{  \pi_S \; : \; S \in S_{\mu}  \}$,
and with every $k$-point $W \in Gr^{\mu}(V) $
we associate the $1$-dimensional $k$-vector subspace
$k \cdot \left\{    \pi_S \mapsto \pi_S(w) \, : \, S \in S_{\mu} \right\}
\subset \Pi (S_{\mu})^* $,
which is the same
for all admissible isomorphisms $w$ of $W$.

Note also that $\pi_S (w) \ne 0 $ if and only if
$W \in Gr^{\mu}_{A_S}(V)$, i.~e. $W \oplus A_S = V$. \\[4pt]

Consider the group
$$GL (V,V_0) \eqdef \{ g \in Aut_k(V) \;:\; \dm\nolimits_k (gV_0 + V_0) /
(gV_0 \cap V_0)   < \infty              \} \quad \mbox{.}$$
Note that "in coordinats"  $GL (V,V_0) $
corresponds to the group of invertible elements
of the algebra of matrices $(A_{i j})_{ i,j \in \sdz}$
such that for every integer $l$ the number of non-zero $A_{i j}$
with $j \ge l$, $i \le l$ is finite.
And the action of $(A_{i j})_{i,j \in \sdz}$ on $k((t))$
is $E_{i j} (t^j) = t^i$. (Here $E_{i j}$ is the matrix with a $1$
on the $(i,j)$-th entry and zeros elsewhere.)
Define the subgroup
$$GL_0 (V, V_0) \eqdef
\{ g \in GL(V,V_0) \; : \;
\dm\nolimits_k gV_0 / (gV_0 \cap V_0) = \dm\nolimits_k V_0 / (gV_0 \cap V_0)
       \}      \quad \mbox{.}
$$
Then we have the obvious action of $GL_0 (V, V_0)$
on $Gr^{\mu}(V)$.
 Let the group $\hat{GL_0} (V, V_0)$
be the central extension of the group $GL_0 (V, V_0)$
by means of $k^*$, which is the group of all
automorphisms of linear bundle $Det$
on $Gr^0(V)$
which cover actions of elements of group $GL_0 (V, V_0)$
on $Gr^0(V)$.
There exists an explicit description of
$\hat{GL_0} (V, V_0)$
and its action on $Det \mid_{Gr^{\mu}(V)}$ (see~\cite{PS}).
Let $H$ be the group
$$
H = \left\{ (g, E) \in GL_0 (V, V_0)  \times GL_k (t^{-1}k[t^{-1}]) \quad :
\quad   \qquad \qquad \qquad \qquad \qquad \qquad \qquad \qquad \qquad \qquad
\right. $$
$$
\left.
\qquad  \qquad  \qquad \qquad \qquad \qquad \qquad \qquad
g_{--} \, - E \;  : \; t^{-1}k[t^{-1}] \to t^{-1}k[t^{-1}]
 \quad \mbox{has  finite rank}
 \right\}  \mbox{,}
$$
where $g_{--}  \eqdef (g)_{i,j \le -1}$.  Let $N \subset H$
be the normal subgroup  defined as $N = \{ (1,E) \in H  \: : \:
\det E = 1               \}$.
Then $\hat{GL_0} (V, V_0) = H/N$,
and the projection on the first factor gives us the map onto
$GL_0 (V, V_0) $.
Now describe the action of $\hat{GL_0} (V, V_0) $
on $Det \mid_{Gr^{\mu}(V)}$. \\[6pt]
\underline{$\mu = 0$} : for  $(g,E) \in H$
and $[\lambda, w] \in Det \mid_{Gr^0 (V)}$
\begin{equation}   \label{f}
(g,E) [\lambda, w]  \eqdef [\lambda, gwE^{-1}]
\end{equation}
This expression gives us the correct action of
$\hat{GL_0} (V, V_0) $
on $Det \mid_{Gr^0 (V)}$. \\[6pt]
\underline{$\mu \ne 0$}.
Define the action $a \in \hat{GL_0} (V, V_0) $
on $Det \mid_{Gr^{\mu} (V)}$
as action of
\begin{equation}  \label{ff}
\sigma^{- \mu} \cdot \tilde{\sigma}^{\mu} (a) \cdot \sigma^{\mu}
\quad \mbox{,}
\end{equation}
where $\sigma: Det \mid_{Gr^{\mu} (V)}  \to Det \mid_{Gr^{\mu-1} (V)} $
is defined by the formula $\sigma \cdot [\lambda, w]  = [\lambda,
t^{-1} \circ w]$; $t^{-1} \circ w \: : \: H_{0(\mu-1)} \to t^{-1}W$,
$t^{-1} \circ w (x) = t^{-1} w (tx)$ is admissible, $x \in H_{0(\mu-1)}$;
and the automorphism $\tilde{\sigma} \: : \:
\hat{GL_0} (V, V_0)  \to \hat{GL_0} (V, V_0) $
is induced by the following endomorphism of $H$:
$$
(g, E) \longmapsto (t^{-1} g t, E_{\sigma})  \quad \mbox{, where}
$$
$$
E_{\sigma} \mid_{t^{-2}k[t^{-1}]} = t^{-1} E t
\qquad
\mbox{and}
\qquad
E_{\sigma}(t^{-1})=1 \quad \mbox{.}
$$ \\[1pt]

Now consider the group $\Gamma = GL(2, k[[z]])$,
which is the subgroup of $GL_0 (V, V_0)$.

\begin{lemma}   \label{zvezd}
The central extension $\hat{GL_0} (V, V_0)$
splits over $\Gamma$.
\end{lemma}
\proof.
Define a group homomorphism
\begin{equation}      \label{fff}
h \; : \; \Gamma \longrightarrow H \quad \mbox{by}  \quad h(\gamma)
= (\gamma ,\gamma_{--}) \quad \mbox{for} \quad \gamma \in \Gamma  \;
\mbox{.}
\end{equation}
(Observe that $\Gamma$ keeps $V_0$ stable and hence
$\gamma_{--} \in GL_{k} (k[t^{-1}]))$.
Then the group homomorphism
	$$
\Gamma \stackrel{h}{\longrightarrow} H \to H/N = \hat{GL_0} (V, V_0)
$$
splits the central extension over $\Gamma$.
\\[10pt]

 In $V = k((z)) \oplus k((z))$ consider  two lines:
$$
l_1= k((z)) \oplus 0   \qquad \mbox{and} \qquad l_2= 0 \oplus k((z))   \quad
\mbox{,}
$$
which correspond to $k((t^2))$ and $t k((t^2))$ in $k((t))$.
For any $S \in S_{\mu}$ let
$$
n_i(S) \eqdef \dm\nolimits_k H_S \cap V_0 \cap l_i -
\dm\nolimits_k l_i / (H_S \cap l_i
+ l_i \cap V_0)   \; \mbox{,} \quad i=1,2
$$
Note that
\begin{equation} \label{twodiez}
n_1(S) + n_2(S) = \mu    \quad \mbox{.}
\end{equation}
Let $
\left(
\begin{array}{cc}
\alpha^{a_1} & 0 \\
0  &   \nu^{a_2}
\end{array}
\right)
\in \Gamma$,
 $\alpha, \: \nu \in k^*$ $a_i \in \dz $.
By lemma~\ref{zvezd} we can consider an action of
$
\left(
\begin{array}{cc}
\alpha^{a_1} & 0 \\
0  &   \nu^{a_2}
\end{array}
\right)
$  on the  determinant bundle $Det$.

\begin{prop}   \label{dddd}
Let $S \in S_{\mu}$, let $\pi_S$ be the corresponding global section of
$Det^{*} \mid_{Gr^{\mu}(V)}$.
Let
$
\left(
\begin{array}{cc}
\alpha^{a_1} & 0 \\
0  &   \nu^{a_2}
\end{array}
\right)^*
$
be the action of
$
\left(
\begin{array}{cc}
\alpha^{a_1} & 0 \\
0  &   \nu^{a_2}
\end{array}
\right)
$
on $H^0 (Gr^{\mu}(V), Det^*)$.
Then
$$
\left(
\begin{array}{cc}
\alpha^{a_1} & 0 \\
0  &   \nu^{a_2}
\end{array}
\right)^*
(\pi_S)=
\alpha^{a_1 \cdot n_1(s)}
\nu^{a_2 \cdot n_2(s)}
\pi_S  \quad \mbox{.}
$$
\end{prop}
\proof follows from
an explicit description of actions of
$
\left(
\begin{array}{cc}
\alpha & 0 \\
0  &  1
\end{array}
\right)
$
and
$
\left(
\begin{array}{cc}
1 & 0 \\
0  &  \nu
\end{array}
\right)
$
on the determinant bundle.
We get this description by means of direct
calculations from formulas~(\ref{f}), (\ref{ff}) and~(\ref{fff}).\\[10pt]

For any $S \in S_{\mu}$ let $n(S)  \eqdef n_1(S) - n_2(S)$.

\begin{lemma}   \label{ch}
Let $W \in Gr^{\mu} (V)$.
Let $S(W) = \{S \in S_{\mu} \; : \; \pi_S(W) \ne 0 \}$.
Then we have
\begin{enumerate}
\item
if $\min\limits_{S \in S(W)} n(S) > -\infty$, then $l_1 \cap W \ne 0$;
\item if $\rank A_W = 1$, $l_1 \cap W \ne 0$,
then $\min\limits_{S \in S(W)} n(S) > - \infty$.
\end{enumerate}
\end{lemma}
\proof.

Remind that
\begin{equation}   \label{kr}
\pi_S (W) \ne 0 \quad \mbox{if and only if}
\quad
W \oplus A_S = V
\end{equation}
Consider the second statement of this lemma.
If $l_1 \cap W  \ne 0$, then from $\dm_k W \cap V_0 < \infty $
we have
$\dm_k W \cap l_1 \cap V_0 < \infty$.
And from $\rank A_W = 1$
we have
$$
\dm\nolimits_k k((z)) / (A_W + k[[z]])  < \infty
\quad \mbox{, therefore}
$$
$$
\dm\nolimits_k l_1 / (W \cap l_1 + V_0 \cap l_1)  < \infty
\quad
\mbox{.}
$$
Now let some $S' \in S(W)$, then
$W \oplus A_{S'} = V$.
Hence $W \cap A_{S'} = 0$. Hence
\begin{equation}   \label{dvezv}
(W \cap l_1) \cap (A_{S'} \cap l) = 0 \quad \mbox{.}
\end{equation}
From~(\ref{dvezv}) we obtain
$$
n_1(S')  \ge \dm\nolimits_k W \cap l_1 \cap V_0 -
\dm\nolimits_k l_1 /( W \cap l_1 + V_0 + l_1) \quad \mbox{.}
$$
From~(\ref{twodiez}) we have
$
n(S') = 2 n_1(S') - \mu
$. Therefore
\begin{equation}  \label{gg}
n(S') \ge 2 \left( \dm\nolimits_k W \cap l_1 \cap V_0
- \dm\nolimits_k l_1 / (W \cap l_1 + V_0 \cap l_1)
 \right) - \mu  \quad  \mbox{.}
\end{equation}
The second statement of this lemma is proved.

Now consider the first statement of this lemma.
From
$\min\limits_{S \in S(W)} n(S) > - \infty$
we see that
$\min\limits_{S \in S(W)} n_1(S) > - \infty$.
Let
$$
\tau = \min(0,\min\limits_{S \in S(W)} n_1(S)) -1   \quad \mbox{.}
$$
For any integer $m>0$ define
$$
B_m  = z^{\tau}k[[z]] \oplus z^{m + \mu} k[[z]] \subset V
\quad \mbox{.}
$$
 Fix some integer $m_0 > \max (0, -\mu, -\tau)$ such that
$$
W \cap (z^{m_0 + \mu} V_0) = 0
\qquad
\mbox{and}
\qquad
V / (W + z^{-m_0} V_0) = 0    \; \mbox{.}
$$
 We will show that for every $m > m_0$
\begin{equation}  \label{zvzv}
W \cap B_m  \ne 0  \; \mbox{.}
\end{equation}
For this goal we consider the images $\tilde{W}$
of $W \cap (z^{-m}V_0)$
and $\tilde{l_1}$ of $l_1 \cap (z^{-m} V_0)$
in the finite-dimensional Grassmanian manifold
$$Gr (2m+\mu, z^{-m}V_0 / z^{m+\mu} V_0)  \; \mbox{.} $$
Now~(\ref{zvzv}) is equivalent to
$$
(\tilde{W} \cap  \tilde{l_1})
\cap
((z^{\tau}V_0 / z^{m+\mu} V_0 ) \cap \tilde{l_1}) \ne 0  \; \mbox{.}
$$
And the last formula follows from  comparison
of dimensions of finite-dimensional vector spaces:
$$
\dm\nolimits_k \tilde{l_1} = 2m + \mu
$$
$$
\dm\nolimits_k (z^{\tau} V_0 / z^{m+\mu} V_0) \cap \tilde{l_1}     =
m+\mu -\tau
$$
\begin{equation} \label{rrr}
\dm\nolimits_k
\tilde{W} \cap \tilde{l_1} >
m + \tau     \quad \mbox{,}
\end{equation}
where~(\ref{rrr}) is true by the following reason. \\
If
\begin{equation} \label{bec}
\dm\nolimits_k \tilde{W} \cap \tilde{l_1} \le m + \tau  \quad \mbox{,}
\end{equation}
then we find a $k$-vector subspace $K \subset \tilde{l_1}$
such that  $K$
is  the span of  $\{ z^{i_k} \}$
for some integers $i_{k}$
and
$\tilde{l_1}  = (\tilde{W} \cap \tilde{l_1}) \oplus K $.
Then we can find $\tilde{L} \in
 Gr(2m+\mu, z^{-m}V_0 / z^{m+\mu} V_0) $
such that  $\tilde{L}$
is the span of  $\{ z^{i_l} \}$ for some integers $i_l$,
$\tilde{L} \cap \tilde{l_1} = K$,
and
$$z^{-m}V_0 / z^{m+\mu}V_0 = \tilde{W} \oplus \tilde{L} \quad \mbox{.}$$

Let $L$  be the pre-image of $\tilde{L}$
in $z^{-m}V_0$.
We can find the Pl\"ucker coordinate
$S' \in S_{\mu}$ such that $L = A_{S'}$. Then,  by construction,
\begin{equation}  \label{krkr}
V = W \oplus A_{S'}       \quad \mbox{,}
\end{equation}
and from~(\ref{bec}) $n_1(S') \le \tau < \min\limits_{S \in S(W)} n_1(S)$.
The last formula contradicts~(\ref{krkr}).
Therefore formula~(\ref{rrr}) is proved. Hence  formula~(\ref{zvzv})
is true.

Now from the Fredholm condition on $W$ we
have $\dm_k W \cap (z^{\tau} V_0)  < \infty$.
From~(\ref{zvzv}) we have non-empty filtration
$ W \cap   B_m $
for all sufficiently large integers m.
The last one is equivalent to
$$
(W \cap (z^{\tau} V_0)) \cap l_1   \ne 0  \quad \mbox{.}
$$
Indeed, assume the converse. Then choose
a finite $k$-basis $e_1, \ldots, e_n$
of $ W \cap z^{\tau}V_0 $
such that $\nu_2(e_i) \ne \nu_2(e_j)$
for all pairs $i,j$. (Here $\nu_2(e_k) = p$ if $e_k \in B_p$
but $e_k \ne B_{p+1}$.)
Then for any $m > \max\limits_i  \nu_2(e_i) $
we have $W \cap B_m = 0$.
This contradiction concludes the proof of  lemma~\ref{ch}. \\[3pt]

\begin{lemma}  \label{chch}
Let $W \in Gr^{\mu}(V)$,
$\rank A_W = 1$,
and condition~1 or condition~2
of the previous lemma is true,
then
\begin{equation} \label{chk}
 \min_{S \in S(W)} n(S) =
2 \left( \dm\nolimits_k W \cap l_1 \cap V_0
 - \dm\nolimits_k l_1/(W \cap l_1 + V_0 \cap l_1)
           \right) -  \mu
		\end{equation}
\end{lemma}
\proof.
From the proof of expression~(\ref{gg})
we have that the left hand side of expression~(\ref{chk})
is greater or equal to the right hand side of~(\ref{chk}).

Now fix some integer $m > 0$ such that
$$
W \cap (z^{m + \mu} V_0) = 0
\qquad
\mbox{and}
\qquad
V / (W + z^{-m} V_0) = 0    \; \mbox{.}
$$
 Consider the images $\tilde{W}$
of $W \cap (z^{-m}V_0)$
and $\tilde{l_1}$ of $l_1 \cap (z^{-m} V_0)$
in
$  z^{-m}V_0 / z^{m+\mu} V_0$.

Then find a $k$-vector subspace $K \subset \tilde{l_1}$
such that  $K$
is  the span of  $\{ z^{i_k} \}$
for some integers $i_{k}$
and
\begin{equation}  \label{l}
\tilde{l_1}  = ( \tilde{W} \cap \tilde{l_1}) \oplus K \quad \mbox{.}
\end{equation}
Then we can find $\tilde{L} \subset
  z^{-m}V_0 / z^{m+\mu} V_0 $
such that  $\tilde{L}$
is the span of $\{ z^{i_l} \}$ for some integers $i_l$,
$\tilde{L} \cap \tilde{l_1} = K$,
and
\begin{equation}   \label{ll}
z^{-m}V_0 / z^{m+\mu}V_0 = \tilde{W} \oplus \tilde{L}  \quad \mbox{.}
\end{equation}
Denote by $L$
 the pre-image of $\tilde{L}$
in $z^{-m} V_0$.
Then there exists some $S' \in S_{\mu}$
such that $L = A_{S'}$  and $V = W \oplus A_{S'}$.
From~(\ref{l}) and~(\ref{ll})
we have $n(S') = n(W)$.   Lemma~\ref{chch}  is proved.

\section{Stable and non-stable points}
For torsion free sheaves  on $C$ one can define the notion
of "stability" and "semistability". (See~\cite{N}.)
Remind these definitions.

\begin{defin}
A rank $2$ torsion free sheaf $\f$ on a curve $C$
is called a semistable sheaf if for any subsheaf $\g \in R(\f)$
$$
2 \chi(\g) \le \chi (\f)   \quad \mbox{.}
$$
\end{defin}
\begin{defin}
A rank $2$ torsion free sheaf $\f$ on a curve $C$
is called a stable sheaf if for any subsheaf $\g \in R(\f)$
$$
2 \chi(\g) < \chi (\f)   \quad \mbox{.}
$$
\end{defin}     \vspace{3pt}

\begin{th}  \label{10vkv}
Let $W = K (C, p, \f, t_p, e_p)$.
Then the sheaf $\f$  is not a (semi)stable sheaf if  and only if
there exists some $g \in SL (2, k[[z]])$
such that for any $S \in S(gW)$  $n(S) \ge 0$  $(>0)$.
\end{th}
\proof follows from proposition~\ref{p1},
lemmas~\ref{Shurik}, \ref{ch} and~\ref{chch}, and  facts that

1) group $SL(2, k[[z]])$ acts transitively on the set of all
$1$-dimensional over $k((z))$ vector subspaces in $V$;

2) the group $SL(2, k[[z]])$ keeps $V_0$  stable.  \\[4pt]

\noindent
{\bf Remark}.
From lemma~\ref{chch} and formula~(\ref{ech}) we can also
obtain numerical datas of nonstability of the sheaf. \\[5pt]

\noindent {\bf Example 1.} (See also~\cite[ prop.~3.8]{M}.)
Let $\mu = 2\nu$ for some integer $\nu$,
and $W = K(C, p, \f, t_p, e_p)$
is a point of the big cell of $Gr^{\mu}_V$, i.~e.,
\begin{equation}  \label{33}
\pi_{0(\mu)} (W) \ne 0    \quad \mbox{,}
\end{equation}
then the sheaf $\f$ is semistable.
Indeed, if condition~(\ref{33}) is true,
then $W \oplus z^{\nu} V_0 = V$
and for any $g \in SL(2, k[[z]])$,
$g (z^{\nu} V_0) = z^{\nu} V_0$,
$gW \oplus z^{\nu} V_0 = V $.
Therefore
$\pi_{0(\mu)} (gW) \ne 0 $   . But $n (0(\mu)) = 0$,
and we apply theorem~\ref{10vkv}. \\[2pt]

 We can generalize example 1  in the following way. \\[1pt]

\noindent {\bf Example 2.}
With every $W \in Gr^{\mu}(V)$
one can associate (see~\cite{P}):

1) the function $T_W : \dz \to \{ 0,1,2 \}$
such that $T_W(i) = 2$ for $i \ll 0 $, $T_W(i)= 0 $ for
 $i \gg 0$, and $T_{gW} = T_W$ for any $g \in SL(2, k[[z]])$
$$
T_W (i) \eqdef \dm\nolimits_k (W \cap z^i V_0)/ (W \cap z^{i+1} V_0) \;
\mbox{,}
$$

2)  the point
\begin{equation}                 \label{resh}
x_W = (x_{W,i_1}, \ldots , x_{W, i_{l_W}})
\in (\dbp^1_k)^{l_W}
   \quad \mbox{,}
\end{equation}
where $i_1, \ldots, i_{l_W}$
are all the integers such that
$ T_W(i_j)=1 $, and the point $x_{W,i_j} \in \dbp^1_k $
is the $1$-dimensional $k$-vector space
$(W \cap z^{i_j} V_0)/(W \cap z^{i_j+1} V_0)$
in the $2$-dimensional $k$-vector space $z^{i_j}V_0 / z^{i_j+1}V_0$.

  Note that all $\dbp^1_k$ from~(\ref{resh})
have the canonical homogeneous coordinates,
which are induced by
$z^{i_j} \oplus 0 \in V$
and   $0 \oplus z^{i_j} \in V$
for the corresponding $j$.
By means of these coordinates we can identify
all $\dbp^1_k$ from~(\ref{resh})
and consider $x_W$ as $l_W$
ordered points on a line.

There exists the natural diagonal action of the group
$SL(2,k)$ on $(\dbp^1_k)^{l_W}$.
Consider the Segre imbedding
$(\dbp^1_k)^{l_W}
\hookrightarrow
\dbp^{2l_W -1}_k
$.
There is a unique linear action of
$SL(2,k)$ on
$\dbp^{2l_W -1}_k
$
which makes the
Segre imbedding into a $SL(2,k)$-morphism
(the $l_W$-fold tensor represantation of $SL(2,k)$).
Also, $SL(2,k)$ is a reductive group.
Therefore G.I.T. determines (semi)stable
points on $(\dbp^1_k)^{l_W}$ with respect to the $SL(2,k)$ action
and the Segre imbedding.
From~\cite{N} we have the following explicit description:
$x \in (\dbp^1_k)^{l_W}$ is a (semi)stable point if and only if
no point of $\dbp^1_k$ occurs as a component of
$x \ge l_W/2$ $(> l_W/2)$ times.

Now we can formulate the following proposition.

\begin{prop} \label{den}
Let  $W = K (C, p, \f, t_p, e_p)$.
Then

1) If $1 \notin \Image T_W$, then $\f$ is a semistable sheaf.

2) If $x_W \in (\dbp^1_k)^{l_W}$
is a (semi)stable point with respect to the $SL(2,k)$ action
and the Segre imbedding, then $\f$ is a (semi)stable sheaf.
\end{prop}
\proof.

 Indeed, if
$1 \notin \Image T_W$,
then let $M \subset V$
be the span of $\Bigl\{ \{z^i \oplus 0   \},  \{0 \oplus z^i  \} \: :  \:
T_W(i) = 2
   \Bigr\}   $.
We can find $S' \in S_{\mu}$
such that $M = H_{S'}$.
By construction, $W \oplus A_{S'} = V$. Hence $\pi_{S'} (W) \ne 0$.
For any $g \in SL(2, k[[z]])$
we have $T_{gW}=T_W$. Therefore
$\pi_{S'} (gW) \ne 0$.
From theorem~\ref{10vkv} and $n(S')=0$
we see that $\f$  is a semistable sheaf.

Now consider the case 2. Let
$
x_W = (x_{W,i_1}, \ldots , x_{W, i_{l_W}})
\in (\dbp^1_k)^{l_W}
$.
For every $k \in \{i_1, \ldots, i_{l_W}  \}$
we can choose $z_k \in V$
equal to $\{  z^k \oplus 0     \}$
or $\{ 0 \oplus z^k  \}$
such that  the image of the $k$-line $k \cdot z_k$
in $z^k V_0 / z^{k+1} V_0$
does not coincide with $x_{W,k}$.
Let  $M \subset V$
be the span of $\Bigl\{ z_k : k \in \{i_1, \ldots, i_{l_W}   \}
       \Bigr\} $
and $\Bigl\{ \{z^i \oplus 0 \},  \{ 0 \oplus z^i \}   \:
: \: T_W(i)=2  \Bigr\}  $.
Then there exists $S' \in S_{\mu}$
such that $M = H_{S'}$
and $W \oplus A_{S'} = V$,
i.~e., $\pi_{S'} (W) \ne 0$.
Moreover, if $x_W \in (\dbp^1_k)^{l_W} $
is a (semi)stable point,
then from the explicit description of
such points we can choose $z_k \in V$
such that
\begin{equation} \label{dom}
n(S') < 0  \; (\le 0)  \quad   \mbox{.}
\end{equation}
The action
of $SL(2,k)$ on   $(\dbp^1_k)^{l_W} $
is induced  by the action of $SL(2,k[[z]])$
on $V$.
Therefore $x_{gW} \in (\dbp^1_k)^{l_{gW}} $
is a (semi)stable point for any $g \in SL(2, k[[z]])$.
Now from this fact, formula~(\ref{dom}),
and theorem~\ref{10vkv}
we get that $\f$
is a (semi)stable sheaf.
\\[4pt]

\noindent
{\bf Remark}.
Suppose that
$x \in
(\dbp^1_k)^{l} $
is not a semistable point with
respect to the $SL(2,k)$ action and the Segre imbedding.
Then it is easy to construct quintets
$(C, p, \f, t_p, e_p)$
with stable, semistable, and nonstable sheaves $\f$
such that for
$W = K(C, p, \f, t_p, e_p)$
we have $l_W = l$ and $x_W = x$.
\\[4pt]

\noindent
{\bf Remark}.
One can introduce the cellular decomposition of $Gr^{\mu}(V)$
indexed by the set $S_{\mu}$ (see~\cite{PS}):
$W \in Gr^{\mu}(V)$
belongs to the cell
$S = \{ s_{-\mu +1}, s_{-\mu +2}, \ldots        \}$
if  $\pi_S (W) \ne 0$
while $\pi_{S'}(W) = 0$
unless $s'_{-\mu +1} \ge s_{-\mu +1},
s'_{-\mu +2 }  \ge s_{-\mu +2}, \ldots $.
(Such $S \in S_{\mu}$
is uniquely defined for each $W \in Gr^{\mu}(V)$.)
It is not difficult to see that the case~1 of proposition~\ref{den}
corresponds to the following statement:
if $W = K (C, p, \f, t_p, e_p)$
belongs to the cell $S \in S_{\mu}$
such that $H_S \cap l_1 = H_S \cap l_2$,
then $\f$ is a semistable sheaf.\\[5pt]

 Now note that the action of the group $SL(2, k[[z]])$
 on $Gr^{\mu}(V)$ keeps
$W = K(C, p, \f, t_p, e_p) \in Gr^{\mu} (V)$
in the image of the Krichever map and does not change
the curve $C$, the point $p$,
the local parameter $t_p$, and the sheaf $\f$.
So this action is an action
 only on "basises $e_p$".
 Therefore if we want to
 get something like "moduli spaces"
 of sheafs on $C$,
 we would need "to kill"
 the action of $SL(2, k[[z]])$
 on $Gr^{\mu} (V)$.

 On the other hand, by lemma~\ref{zvezd},
 the group  $SL(2, k[[z]])$
 acts on $Det \mid_{Gr^{\mu} (V)}$.
 Moreover, from this action we obtain
 an
  infinite-dimensional
  representation
  of  $SL(2, k[[z]])$
  in the $k$-vector space $\Pi (S_{\mu})$.
	 (In general, this representation is the restriction of representation of
    $\hat{GL_0} (V, V_0) $
     in $\Pi (S_{\mu})$, which can be got in the same way and is usually
  called the infinite wedge-representation (see~\cite{KP}).)
  Hence we have a linear action of
  $SL(2, k[[z]])$    on $\dbp (\Pi (S_{\mu})^*)$.
  And by construction, this action
  is compatible with the action of
  $SL(2, k[[z]])$
  on $Gr^{\mu} (V)$
  via the Pl\"ucker imbedding.
  These arguments lead us to the following definitions.

  For any $W \in Gr^{\mu} (V)$
  denote by $\hat{W}$
  some nonequal to $0$ element from   $ \Pi (S_{\mu})^*  $
  such that the image of $\hat{W}$
  in $\dbp (\Pi (S_{\mu})^*)$
  coincides with the image of $W$
  via the Pl\"ucker imbedding.

  Let
	$$
	T = \left\{
  \left(
\begin{array}{cc}
\lambda  & 0 \\
0  &   \lambda^{-1}
\end{array}
\right)
    \in
  SL(2,k[[z]]),
  \: \lambda \in k^*
	       \right\}
	$$
  be a 1-parametric subgroup:
  $
  k^*  \to SL(2, k[[z]])
  $.
  For any $g \in SL(2, k[[z]] $
  define the group  ${}^gT = gTg^{-1}$.
  Let ${}^gT(\lambda) = g
  \left(
\begin{array}{cc}
\lambda  & 0 \\
0  &   \lambda^{-1}
\end{array}
\right)
	g^{-1}$  be an element from ${}^gT$.

  It will be convinient
 to identify $k^*$
 with a subset of
 $\dbp^1_k$ in the obvious manner.
 In fact, for any $\lambda \in k$,
 we identify  $\lambda$ with the point $(1, \lambda)$
 in $\dbp^1_k$,
 and write $\infty$ for the extra point $(0,1)$.
 Note that the morphism $k^* \to \Pi (S_{\mu})^*  \; : \;
 \lambda \mapsto {}^gT(\lambda) \hat{W}$
 may or may not extend to the points $0$, $\infty$,
 but if it does extend to either or both of these
 points the extension is unique.
 Thus we can define  the expressions
  $
  \Lim\limits_{\lambda \to 0}
  {}^gT(\lambda) \hat{W}
  $
	  and
  $
  \Lim\limits_{\lambda \to \infty}
  {}^gT(\lambda) \hat{W}
  $ in an obvious way.

  \begin{defin}
  A closed point
  $W \in Gr^{\mu} (V)$
  is a semistable point with respect
  to the 1-parametric subgroup ${}^gT$
  if  for $x = 0$ and $x = \infty $
	$
  \lim\limits_{\lambda \to x}
  {}^gT(\lambda) \hat{W}
  $
	 does not exist
  or
  $
  \Lim\limits_{\lambda \to x}
  {}^gT(\lambda) \hat{W}   \ne 0
  $.
   \end{defin}
  \begin{defin}
  A closed point
  $W \in Gr^{\mu} (V)$
  is a stable point with respect
  to the 1-parametric subgroup ${}^gT$
  if  for $x = 0$ and $x = \infty $
	$
  \Lim\limits_{\lambda \to x}
  {}^gT(\lambda) \hat{W}
  $
	 does not exist.
   \end{defin}
  Now using these definitions and  proposition~\ref{dddd}
  we can reformulate  theorem~\ref{10vkv}
  in the following way.(Compare~\cite{N}.)
  \begin{th}       \label{th1}
   Let $(C,p,\f, t_p, e_p)$
   be a quintet,
   $\chi (\f) =  \mu$.
   Then the sheaf $\f$
   is a (semi)stable sheaf on the curve $C$
   if and only if the point
   $K(C,p,\f, t_p, e_p) \in Gr^{\mu}(V)$
   is a (semi)stable point
   with respect to the 1-parametric
   subgroups ${}^gT$ for all $g\in SL(2, k[[z]])$.
  \end{th}

\noindent {\bf Acknowledgments}  \\
 This note was done during my visit to the Abdus Salam ICTP
in the end of 1998. I am very much grateful to Professor M.~S.~Narasimhan
for the invitation me to the ICTP and for the excellent working
conditions.
I am deeply grateful also to my scientific advisor  Professor
A.~N.~Parshin for the permanent attention and help.


\vspace{5pt}

{\em e-mail: } $d_{-}osipov@chat.ru$
\end{document}